# New RBF collocation methods and kernel RBF with applications

W. Chen

Department of Informatics, University of Oslo, P.O.Box 1080, Blindern, 0316 Oslo, Norway

**Abstract.** A few novel radial basis function (RBF) discretization schemes for partial differential equations are developed in this study. For boundary-type methods, we derive the indirect and direct symmetric boundary knot methods. Based on the multiple reciprocity principle, the boundary particle method is introduced for general inhomogeneous problems without using inner nodes. For domain-type schemes, by using the Green integral we develop a novel Hermite RBF scheme called the modified Kansa method, which significantly reduces calculation errors at close-to-boundary nodes. To avoid Gibbs phenomenon, we present the least square RBF collocation scheme. Finally, five types of the kernel RBF are also briefly presented.

## 1. Introduction

Many existing meshfree methods require using the moving least square (MLS). Exceptionally, the numerical schemes based on the radial basis function (RBF) do not need the MLS at all and are inherently meshfree. For some new advances on the RBF see Buhmann's excellent survey [1]. The RBF has physical backgrounds of field and potential theory [2] and is justified mathematically by integral equation theory [3]. Among RBF numerical schemes, famous are the Kansa method [4], Hermite symmetric RBF collocation method [5,6] and the method of fundamental solution (MFS) [7].

   The Kansa's method is the very first domain-type RBF collocation scheme with easy-to-use merit, but the method lacks symmetric interpolation matrix due to boundary collocation. The Hermite RBF collocation method kills the unsymmetrical drawback. Like the Kansa's method, however, the method suffers relatively lower accuracy in boundary-adjacent region. The MFS, also known as the regular BEM, is a simple and efficient boundary-type RBF scheme, but the controversial artificial boundary outside physical domain hinders its practical applications. The boundary knot method, recently introduced by the present author [2,3], surpasses the MFS in that it employs the nonsingular general solution instead of the singular fundamental solution and thus no longer requires the arbitrary fictitious boundary. Albeit better than the MKM, the BKM loses symmetric merit whenever the presence of mixed boundary conditions. It is worth pointing out that all these RBF schemes are indirect and global. The indirect methods mean that the expansion coefficients rather than physical variables are used as the basic variable, while the global interpolation causes the ill-conditioning interpolation matrix and susceptible to Gibbs phenomenon amid weak continuity of physical solution.
   The purpose of this paper is to introduce a few new RBF discretization schemes of boundary and domain types to overcome the aforementioned shortcomings. On the other hand, the proper RBF is also an essential issue leading to an efficient and stable solution. In general, there is not an operational approach to create efficient RBF available now. Based on the underlying relationship between the RBF and the Green integral [3], we summarize five approaches constructing kernel RBF [8], which also cover all existing popular RBFs.

The rest of this paper is structured as follows. In section 2, we establish the symmetric BKM and the direct BKM, and then, the symmetric boundary particle method (BPM) is developed by using the multiple reciprocity principle [9]. Unlike the BKM, the BPM does not require the interior nodes for inhomogeneous problems. By using Green integral, section 3 presents the symmetric modified Kansa's method (MKM) to significantly improve the solution accuracy at nodes neighboring boundary. Furthermore, we briefly describe the spline version of the MKM called as the finite knot method (FKM) to produce the sparse symmetric interpolation matrix. Section 4 is concerned with the least square RBF collocation method which may better approximate the solution of lower continuity to avoid Gibbs phenomenon. In section 5 we discuss the kernel RBF, which is the best choice in the approximate expression of particular solution and the domain-type RBF methods. The first-order regulation condition is also proposed for numerical PDE. Finally, section 6 indicates references on numerical experiments and gives reasoning of naming brand-new methods and kernel RBF.

## 2. Boundary-type RBF schemes

### 2.1. Symmetric boundary knot methods

The aim of this section is to derive the symmetric Hermite BKM and direct BKM [9,10]. The following example serves as an illustrative example:

$$L\{u\} = f(x), \qquad x \in \Omega, \qquad (1)$$
$$u(x) = R(x), \qquad x \subset S_u, \qquad (2)$$
$$\frac{\partial u(x)}{\partial n} = N(x), \qquad x \subset S_T, \qquad (3)$$

where $x$ means multi-dimensional independent variable, and $n$ is the unit outward normal. The solution of Eq. (1) can be expressed as

$$u = u_h + u_p, \qquad (4)$$

where $u_h$ and $u_p$ are the homogeneous and particular solutions, respectively. The latter satisfies

$$L\{u_p\} = f(x) \qquad (5)$$

but does not necessarily satisfy boundary conditions. To evaluate the particular solution, the inhomogeneous term is approximated first by

$$f(x) \cong \sum_{j=1}^{N+L} \alpha_j \varphi(r_j), \qquad (6)$$

where $\alpha_j$ are the unknown coefficients. $N$ and $L$ are respectively the numbers of knots on the domain and boundary. The use of interior points here is usually necessary to guarantee the accuracy and convergence of the BKM solution. $r_j = \|x - x_j\|$ represents the Euclidean distance norm, and $\varphi$ is the radial basis function. By forcing approximation representation (6) to exactly satisfy Eq. (5) at all nodes, we can evaluate α. Finally, we can get particular solutions at any point by summing localized approximate particular solutions

$$u_p = \sum_{j=1}^{N+L} \alpha_j \phi(\|x - x_j\|). \qquad (7)$$

On the other hand, the homogeneous solution $u_h$ has to satisfy both governing equation and boundary conditions. Unlike the dual reciprocity BEM (DR-BEM) [11] and MFS [7] using the singular fundamental solution, the BKM [2,3] approximates homogeneous solution by means of nonsingular general solution

$$u_h(x) = \sum_{k=1}^{L} \lambda_k u^{\#}(r_k), \qquad (8)$$

where $k$ denotes index of source points on boundary; $u^{\#}$ represents the nonsingular general solution of operator $L$, and $\lambda_k$ are the desired coefficients. The non-singular general solutions of some frequently-used operators are listed in [2,3,10,12]. The representation (8) will lead to an unsymmetric BKM formulation. By analogy with the symmetric Hermite RBF collocation method presented in [5,6], we modify the BKM approximate expression (8) as

$$u_h(x) = \sum_{s=1}^{L_D} \lambda_s u^{\#}(r_s) - \sum_{s=L_D+1}^{L_D+L_N} \lambda_s \frac{\partial u^{\#}(r_s)}{\partial n}, \qquad (9)$$

where $n$ is the unit outward normal as in boundary condition (3), and $L_D$ and $L_N$ are respectively the numbers of knots at the Dirichlet and Neumann boundary surfaces. The minus sign associated with the second term is due to the fact that the Neumann condition of the first order derivative is not self-adjoint. Collocating at all boundary and interior knots in terms of representation (9), we have

$$\sum_{s=1}^{L_D} \lambda_s u^\#(r_{is}) - \sum_{s=L_D+1}^{L_D+L_N} \lambda_s \frac{\partial u^\#(r_{is})}{\partial n} = R(x_i) - u_p(x_i), \tag{10}$$

$$\sum_{s=1}^{L_D} \lambda_s \frac{\partial u^\#(r_{js})}{\partial n} - \sum_{s=L_D+1}^{L_D+L_N} \lambda_s \frac{\partial^2 u^\#(r_{js})}{\partial n^2} = N(x_j) - \frac{\partial u_p(x_j)}{\partial n}, \tag{11}$$

$$\sum_{s=1}^{L_D} \lambda_s u^\#(r_{ls}) - \sum_{s=L_D+1}^{L_D+L_N} \lambda_s \frac{\partial u^\#(r_{ls})}{\partial n} = u_l - u_p(x_l), \tag{12}$$

where *i, j,* and *l* indicate response knots respectively located on boundary $S_u$, $S_\Gamma$, and domain $\Omega$. The solution of the above simultaneous equations can be decomposed into two steps. The first is to evaluate the unknown boundary expansion coefficients $\lambda$ by using symmetric equations (10) and (11), and then the interior node solution $u_l$ is calculated by Eq. (12). We can employ the obtained expansion coefficients $\lambda_k$ and inner knot solutions $u_l$ to calculate the BKM solution at any knot. It is stressed here that the MFS could not produce the symmetric interpolation matrix in any way.

The above BKM uses expansion coefficients rather than the direct physical variable in the approximation of boundary value. Such BKM is called the indirect BKM. Note that the Neumann condition $N(x)$ at $x \subset S_u$ and Dirichlet condition $D(x)$ at $x \subset S_T$ are unknown in contrast to the prescribed boundary condition (2,3). To simplify the presentation, we use the $D_u$, $N_u$ and $D_\Gamma$, $N_\Gamma$ respectively represent the Dirichlet and Neumann values at $x \subset S_u$ and $x \subset S_T$. In terms of Eqs. (10,11,12), we have

$$A\lambda = \begin{Bmatrix} D_u \\ N_\Gamma \end{Bmatrix}, \qquad B\lambda = \begin{Bmatrix} N_u \\ D_\Gamma \end{Bmatrix}. \tag{13}$$

So

$$\begin{Bmatrix} N_u \\ D_\Gamma \end{Bmatrix} = BA^{-1} \begin{Bmatrix} D_u \\ N_\Gamma \end{Bmatrix}. \tag{14}$$

Note that matrix *A* are symmetric and no indirect expansion coefficients are involved in Eq. (14). The solution procedure is a symmetric direct BKM strategy. Since the general solutions of Helmholtz problem tend to zero at infinity, the BKM is applicable to exterior unbounded Helmholtz problems.

Like the comparisons between the direct and indirect BEMs, the direct BKM has the advantages to problems with sharp corners since the fictitious expansion coefficients may tend to infinity as nodes increase, even if the physical quantity remains well behaved. If an iterative technique is applied, it is also much easier to find a nice initial guess solution for the discretization systems of the direct method than that of the indirect method.

## 2.2. Boundary particle methods

The multiple reciprocity BEM (MR-BEM) applies the multiple reciprocity principle to circumvent the domain integral without using any inner nodes [13]. The shortcoming is uneasily applied to nonlinear problems and requires relatively higher computing effort. In this section, we develop a truly boundary-only RBF scheme based on the multiple reciprocity principle [10].

The MRM assumes that the particular solution of Eq. (1) can be approximated by higher-order homogeneous solution, namely,

$$u = u_h^0 + u_p^0 = u_h^0 + \sum_{m=1}^{\infty} u_h^m, \qquad (15)$$

where superscript $m$ is the order index of homogeneous solution. Through an incremental differentiation operation via operator $\Re\{\}$, we have successively higher order differential equations:

$$\begin{cases} u_h^0(x_i) = R(x_i) - u_p^0(x_i) \\ \dfrac{\partial u_h^0(x_j)}{\partial n} = N(x_j) - \dfrac{\partial u_p^0(x_j)}{\partial n} \end{cases}, \qquad (16)$$

$$\begin{cases} \Re^{n-1}\{u_h^n(x_i)\} = \Re^{n-2}\{f(x_i)\} - \Re^{n-1}\{u_p^n(x_i)\} \\ \dfrac{\partial \Re^{n-1}\{u_h^n(x_j)\}}{\partial n} = \dfrac{\partial (\Re^{n-2}\{f(x_j)\} - \Re^{n-1}\{u_p^n(x_j)\})}{\partial n} \end{cases}, \qquad n=1,2,\ldots, \quad (17)$$

where $\Re^n\{\}$ denotes the $n$-th order of operator $\Re\{\}$, say $\Re^0\{\}=\Re\{\}$ and $\Re^1\{\}=\Re\Re^0\{\}$, $i$ and $j$ are respectively Dirichlet and Neumann boundary knots. $u_p^n$ is the $n$-th order of particular solution approximated by

$$u_p^n = \sum_{m=n+1}^{\infty} u_h^m. \qquad (18)$$

The m-order homogeneous solution is represented by a Hermite expansion

$$u_h^m(x) = \sum_{s=1}^{L_D} \beta_s u_m^\#(r_s) - \sum_{s=L_D+1}^{L_D+L_N} \beta_s \dfrac{\partial u_m^\#(r_s)}{\partial n}, \qquad (19)$$

where $u_m^\#$ is the corresponding $m$-th order fundamental or general solutions. Collocating boundary equations (16, 17), we have discretization equations

$$\sum_{s=1}^{L_D} \beta_s u_0^{\#}(r_{is}) - \sum_{s=L_D+1}^{L_D+L_N} \beta_s \frac{\partial u_0^{\#}(r_{is})}{\partial n} = R(x_i) - u_p^0(x_i)$$

$$\sum_{s=1}^{L_D} \beta_s \frac{\partial u_0^{\#}(r_{js})}{\partial n} - \sum_{s=L_D+1}^{L_D+L_N} \beta_s \frac{\partial^2 u_0^{\#}(r_{js})}{\partial n^2} = N(x_j) - \frac{\partial u_p^0(x_j)}{\partial n} \quad (20)$$

$$\sum_{s=1}^{L_D} \beta_s \Re^{n-1}\{u_n^{\#}(r_{is})\} - \sum_{s=L_D+1}^{L_D+L_N} \beta_s \frac{\partial \Re^{n-1}\{u_n^{\#}(r_{is})\}}{\partial n} = \Re^{n-2}\{f(x_i)\} - \Re^{n-1}\{u_p^n(x_i)\}$$

$$\sum_{s=1}^{L_D} \beta_s \frac{\partial \Re^{n-1}\{u_n^{\#}(r_{js})\}}{\partial n} - \sum_{s=L_D+1}^{L_D+L_N} \beta_s \frac{\partial^2 \Re^{n-1}\{u_n^{\#}(r_{js})\}}{\partial n^2} = \frac{\partial(\Re^{n-2}\{f(x_j)\} - \Re^{n-1}\{u_p^n(x_j)\})}{\partial n}$$

$$n=1,2,\ldots, \quad (21)$$

In terms of the MRM, the successive process is truncated at some order $M$. The practical solution procedure is a reversal recursive process:

$$\beta_k^M \to \beta_k^{M-1} \to \cdots \to \beta_k^0. \quad (22)$$

It is noted that due to

$$\Re^{n-1}\{u_h^n(r_k)\} = u_h^0(r_k), \quad (23)$$

the coefficient matrices of all successive equation are thus the same $Q$

$$Q\beta_k^n = b^n, \ n=M,M\text{-}1,\ldots,1,0, \quad (24)$$

Thus, the LU decomposition algorithm is very suitable for this task. Finally, the solution at any node is given by

$$u(x_p) = \sum_{n=0}^{M} \left( \sum_{s=1}^{L_D} \beta_s u_m^{\#}(r_{ps}) - \sum_{s=L_D+1}^{L_D+L_N} \beta_s \frac{\partial u_m^{\#}(r_{ps})}{\partial n} \right). \quad (25)$$

The BPM can use either singular fundamental solution or nonsingular general solution, respectively relative to the MFS and BKM. The only difference between the BKM (MFS) and BPM lies in how to evaluate the particular solution. The former applies the dual reciprocity principle, while the latter employs the multiple reciprocity principle. The BPM with $M=1$ degenerates into the BKM or MFS without using the inner nodes. The advantage of the BPM over the BKM is that it dose not require interior nodes which may be especially attractive in such problems as moving boundary, inverse problems, and exterior problems. However, the BPM may be more mathematically complicated and computationally costly due to the iterative use of higher-order fundamental or general solutions. It is expected that like the MR-BEM, the truncated order $M$ in the BPM may not be large (usually two or three orders) in a variety of practical uses. Following strategy for direct BKM it is very straightforward to develop the direct BPM methodology. The BKM, BPM and BEM alike always produce full matrix, but the former two are more attractive than the BEM in terms of accuracy, simplicity and efficiency.

## 3. Modified Kansa method based on the second Green identity, spline approximation FKM and direct MKM

Despite great effort, the rigorous mathematical proof of the solvability of the Kansa's method is still missing [14]. The boundary conditions also destroy the symmetricity of its interpolation matrix. As an alternative, refs. 5 and 6 present the symmetric Hermite RBF collocation scheme with sound mathematical analysis of solvability. One common issue in the Kansa's method and symmetric Hermite method, however, is that the numerical solutions at nodes adjacent to boundary deteriorate (by one to two orders) compared with those in central region. Fedoseye et al. [15] propose the PDE collocation on the boundary (PDECB) to effectively remove this shortcoming. The strategy requires an additional set of nodes (inside or outside of the domain) adjacent to the boundary. Like the fictitious boundary of the MFS, the arbitrary placing of these additional nodes may give rise to some troublesome issues. The PDECB also lacks explicit theoretical endorsement. In fact, a similar strategy has been independently proposed by Zhang et al. [16], which collocates both governing and boundary equations on the same boundary nodes. However, the method is unsymmetrical and still lacks explicit theoretical basis.

By using the Green second identity, the solution of Eqs. (1,2,3) is given by

$$u(x) = \int_\Omega f(z) u^*(x,z) d\Omega + \int_\Gamma \left\{ u \frac{\partial u^*(x,z)}{\partial n} - \frac{\partial u}{\partial n} u^*(x,z) \right\} d\Gamma, \quad (26)$$

where $u^*$ is the fundamental solution of differential operator $L\{\}$. $z$ denotes source point. It is noted that the first and second terms of Eq. (26) are respectively equivalent to the particular and homogeneous solutions. If a numerical integral scheme is used to analogize Eq. (26), we have

$$u(x) \cong \sum_{k=1}^{N+L} \omega(x, x_k) f(x_k) u^* + \sum_{k=N+1}^{N+L} Q(x, x_k) \left[ u \frac{\partial u^*}{\partial n} - \frac{\partial u}{\partial n} u^* \right], \quad (27)$$

where $\omega(x, x_j)$ and $Q(x, x_j)$ are the integration weight functions dependent on the integral schemes. Perceiving the RBF as an approximate Green function, we can construct the following interpolation formula

$$u(x) = \sum_{k=1}^{N+L} \alpha_k L^*\{\varphi(r_k)\} + \sum_{k=1}^{N+L_D} \beta_k \varphi(r_k) + \sum_{k=N+L_D+1}^{N+L} \beta_k \left( -\frac{\partial \varphi(r_k)}{\partial n} \right), \quad (28)$$

where $L^*\{\}$ reverses some signs of odd-order derivatives in $L\{\}$ if the latter is not self-adjoint. Note that the boundary nodes are here interpolated twice. The above RBF interpolation scheme differs from the PDECB in that it does not require auxiliary nodes at all and is derived naturally from the Green second identity. Therefore, theoretical and operational ambiguities in the PDECB are eliminated at all. Contrast to the method given in [16], collocating Eqs. (1,2,3) via representation (28) leads to theoretically solid and symmetric modified Kansa's method formulation.

The MKM holds symmetric property with radically symmetric RBF $\varphi$ even if operator $L\{\}$ is not self-adjoint, where we need to create the corresponding operator $L^*\{\}$ with a different sign before the odd-order derivative. For varying parameter problems, the MKM also keeps symmetric merit. For more related details see ref. 9.

The MKM and preceding symmetric BKM establish an underlying connection due to the fact that both are based on the second Green identity. The MKM may be less accurate than the BKM since the latter employs the analytical general solution to approximate homogeneous solution. However, the MKM can produce the symmetric interpolation matrix no matter whether or not the right-hand inhomogeneous term in Eq. (1) includes the dependent variable $u$ and operator $L\{\}$ is self-adjoint.

One can note that the MKM applies the boundary conditions in a similar fashion of the functional FEM. Accordingly the basic strategy localizing the MKM is to apply spline RBF interpolation and domain decomposition, interpreted as influential domain, namely, the solution at one node within a specified subregion is approximated only by those nodes of the same subregion. The subregion can by no way been understood as grids in the standard FEM. The method is still truly meshless. In addition, the strategy does not involve the overlapping across different subregions. We call this localizing MKM as the finite knot method, which produces a sparse banded system matrix. For more details on the FKM and direct MKM see ref. 9.

## 4. Least square RBF collocation method

For problems with discontinuous or weekly continuous solutions, the preceding RBF schemes will encounter the Gibbs phenomenon, which spoils accuracy and stability. It is well known that the least square approach outperforms well the interpolation method in this regard.

Still consider the case of Eqs. (1,2,3). If we choose $N$ source and $M$ field nodes across whole computational domain, where $N$ and $M$ are not necessarily equal and source and field nodes are also not necessarily coincidental, we have the RBF approximate expression

$$u(x) = \sum_{k=1}^{N} \beta_k \psi(r_k), \qquad (29)$$

where $\psi$ is the RBF with the Kansa' method, Hermite RBF method, BKM, BPM, MKM or FKM. Collocating via expression (29) in $M$ field nodes yields

$$G_{M\times N}\beta_{N\times 1} = b_{M\times 1}, \qquad (30)$$

where $G$ is coefficient matrix. Eq. (30) can be solved by the least square approach. The $L_2$ norm of residual errors is

$$\sigma = \sum_{i=1}^{M}\left(\sum_{k=1}^{N}g_{ik}\beta_k - b_i\right)^2, \qquad (31)$$

where $g_{ik}$ are entries of matrix $G$. By forcing the first order derivative of $\sigma$ with respect to $\beta_k$ zero, we have the least square formulation:

$$\hat{G}_{N\times N}\beta_{N\times 1} = \hat{b}_{N\times 1}. \qquad (32)$$

We call this technique the least square RBF collocation method. The methodology is less sensitive to the discontinuity of the physical solution. $L_1$ (Chebyshev) norm of residual errors can also be employed to design a least square RBF scheme. The extension of the methodology to the Galerkin-type method is straightforward. As alternative of RBF interpolation, the application of least square RBF to multivariate data processing is also promising.

## 5. Kernel RBF

Chen [2,3] presented the kernel RBF-creating strategy based on the second Green identity, i.e., formulas (26,27). Chen [8] further proposed three types of kernel RBF. The first is to apply $r^{2m}$ augmented term to enhance the smoothness and ensures sufficient degree of differential continuity since the fundamental solution has a singularity at origin. The TPS is a notable example in this regard. The second strategy is simply the higher order of fundamental and general solutions. The third approach is to replace distance variable $r$ in fundamental or general solutions by $\sqrt{r^2 + c^2}$, where $c$ is shape parameter.

Next is a physical explanation of shape parameter. All complete fundamental solutions consist of essential and complementary elementary functions [17]. The standard singular fundamental solutions used in the BEM involve only the essential part. The complementary terms of the complete fundamental solutions are often understood the nonsingular general solution in terms of the BKM [2,3]. The shape parameter $c$ can be interpreted as the scaling parameter in the simplified form of the complete fundamental solutions and leads to infinite smoothness at the cost of one complementary term. For instance, the MQ and reciprocal MQ are respectively related with general fundamental solutions of 1D and 3D Laplacian. The tricky choose of the shape parameter coincides with skillful implementation of general fundamental solutions. For some details on the kernel RBF see ref. [8].

Following the basic idea of the corrected reproducing kernel approximation [18], next is involved a compactly-supported kernel RBF (CSK-RBF). Multidimensional function $u$ is represented by

$$\bar{u}(\bar{x}) = \int_{\Omega_\xi} w_d(\bar{x} - \bar{\xi}) u(\bar{\xi}) d\Omega_\xi , \qquad (33)$$

where $w_d$ is the corrected kernel function

$$w_c(\bar{x} - \bar{\xi}) = C(\bar{x}, \bar{\xi}) w(\bar{x} - \bar{\xi}). \qquad (34)$$

The correction function $C$ can be a polynomial expansion just like those in [18] or the RBF. For instance, the MQ can be used as correction function to decide local optimal shape parameter by establishing the reproducing conditions (moment conditions). Instead of using polynomial expansion, the point of CSK-RBF is to choose the kernel RBF as the kernel function $w$. Note that $\bar{x} - \bar{\xi}$ implicates not only the Euclidean distance but also wider translation invariance likelihood. The choice of kernel RBF depends on the property of function $u$ or PDE type. The CSK-RBF is seen as a generalized Green integral representation with $u$ itself instead of $f(x)$ in Eq. (26) while keeping higher-order local consistency via correct function. On the other hand, some kernel RBF itself is compactly-supported. For example, those of modified Helmholtz operator. Chen [10,12] has also recently found the higher-order fundamental and general solutions of convection-diffusion, biharmonic, Winkler plate and Burger plate equations.

The kernel RBF also includes time-space RBF [3] using transient fundamental solution and general solution. d'Alembert RBF general solution of 3D linear wave equation is

$$u(r,t) = \frac{1}{r}[f(r+ct) + g(r-ct)]. \qquad (35)$$

Accordingly we construct the characteristic RBF $\varphi(r\text{-}ct)$ for hyperbolic wave problems, where we can choose Laplacian kernel RBF as $\varphi$.

The RBF is closely related to the solution of a PDE with spherically symmetric unbounded domain without boundary conditions. In that sense, the first-order regulation condition means that only the first order derivative of the RBF at origin does not tend to infinity. The higher differential continuity is not necessary for higher-order PDE but may increase the accuracy in the case of smooth solution.

Partial differential equation is the very basic language describing the universe. So it will be a natural choice making the RBF underlying PDE. The operator-dependent kernel RBF is strongly recommended in the solution of PDE and data processing. For example, $e^{-\alpha r}$ is much better than Guassian for diffusion and convection-diffusion problems due to its underlying physical grounds. The numerical integration of the boundary integral equation theory (distribution theory) may be a powerful tool to mathematical analysis of the kernel RBF.

## 6. Numerical experiments and naming of brand-new methods

Due to very limited space, the readers are advised to find the numerical experiments applying the BKM, BPM and MKM with kernel RBF in [10,12,19]. A more complete description of the methods and references go to ref. 9.

The BKM and BPM are two global RBF collocation techniques using nonsingular general solution corresponding to the DR-BEM and MR-BEM. Thus, they may be renamed as the dual reciprocity method of general solution (DR-MGS), multiple reciprocity method of general solution (MR-MGS), and multiple reciprocity MFS(MG-MFS). If combined with particular solution method, we could have PS-MGS just like PS-MFS. Contrast to the known Kansa's method, the domain-type modified Kansa's method borrows the idea of the BKM that the solution of PDE is a sum of particular and homogeneous solutions. A combination of the RBF and least square technique is unknown to the author so far. The least square RBF collocation method reflects the essential components of the technique. The kernel RBF is due to the fact that it applies kernel function of integral equation, especially Green integral, as essential element.